\documentclass[12pt]{amsart}
\usepackage[english]{babel}
\usepackage{egothic}
\usepackage[T1]{fontenc}

\usepackage{amsmath}
\usepackage{amsthm}
\usepackage{amssymb}
\usepackage{amsfonts}
\usepackage[all]{xy} 
\usepackage{mathrsfs}
\usepackage{enumerate}
\usepackage{physics}
\usepackage[notcite, final, notref]{showkeys}
\usepackage{dsfont}
\usepackage[dvips]{color}
\newtheorem{theorem}{Theorem}[section]

\newtheorem{lemma}[theorem]{Lemma}
\newtheorem{proposition}[theorem]{Proposition}

\newtheorem{definition}[theorem]{Definition}

\theoremstyle{definition}
\newtheorem{remark}[theorem]{Remark}

\newtheorem{example}[theorem]{Example}

\newcommand{\BC}{{\mathbb B}{\mathbb C}}
\newcommand{\D}{{\mathbb D}}

\renewcommand{\i}{{\bf i}}
\renewcommand{\j}{{\bf j}}
\renewcommand{\k}{{\bf k}}
\newcommand{\C}{{\mathbb C}}
\newcommand{\R}{{\mathbb R}}
\newcommand{\e}{{\bf e}}
\newcommand{\edag}{{\bf e^\dagger}}

\renewcommand{\Re}{\mathrm{Re}}

\usepackage{xcolor}

\begin{document}

\title[Bicomplex-Real (BC-R) Calculus and It\^o-Hermite]{The Bicomplex-Real Calculus and Applications to Bicomplex Hermite-It\^o Polynomials}

\author[D. Alpay]{Daniel Alpay}
\address{(DA) Schmid College of Science and Technology \\
Chapman University\\
One University Drive
Orange, California 92866\\
USA}
\email{alpay@chapman.edu}

\author[K. Diki]{Kamal Diki}
\address{(KD) Schmid College of Science and Technology \\
Chapman University\\
One University Drive
Orange, California 92866\\
USA}
\email{diki@chapman.edu}

\author[M. Vajiac]{Mihaela Vajiac}
\address{(MV) Schmid College of Science and Technology \\
Chapman University\\
One University Drive
Orange, California 92866\\
USA}
\email{mbvajiac@chapman.edu}

\keywords{ Bicomplex numbers, C-R calculus, B-C-R calculus, Bicomplex gradient operator, It\^o-Hermite polynomials}%
\subjclass{Primary 30G35; Secondary 33C45}%

\thanks{D. Alpay, K. Diki and M. Vajiac thank Chapman University for the Faculty Opportunity Fund which helped fund this research.\ D. Alpay also thanks the Foster G. and Mary McGaw Professorship in Mathematical Sciences, which supported his research. \ K. Diki also thanks the Grand Challenges Initiative (GCI) at Chapman University, which  supported his research. }

%\tableofcontents
\begin{abstract}
In this paper we extend notions of complex C-R-calculus and complex Hermite polynomials to the bicomplex setting and compare the bicomplex polyanalytic function theory to the classical complex case.
\end{abstract}
\maketitle
\section{Introduction}
\setcounter{equation}{0}

In the past decade many applications of hypercomplex analysis and algebra came to light. Among these, for example, the authors have extended the complex perceptron algorithm to the bicomplex case in~\cite{bcplx_perceptron, bcplx_lms}. Inroads in bicomplex function theory and spaces have also been made by many, such as~\cite{alv1, alv2, alss, bcbook, Price}, to cite just a few.
In 1960 Widrow and Hoff(~\cite{WMcB}) have extended the gradient descent technique to the complex domain which was derived with respect to the real and imaginary part. Brandwood in~\cite{Brandwood} used the Wirtinger (or  $\mathbb{C}-\mathbb{R}$)~\cite{Wirtinger} Calculus such that the gradient is considered with respect to complex variables instead. This theory was used to define and use the complex Least Mean Square (LMS) algorithm, first discovered by Widrow and his student Hoff in the real valued case in 1960 (see \cite{WH1960}). 
Later, in 1983, Brandwood~\cite{Brandwood} studied properties of the complex gradient operator using \textit{"Complex-Real-Calculus"} (i.e.  C-R-Calculus or Wirtinger calculus, see also \cite{Kreutz}). There, in~\cite{Brandwood}, Brandwood applied these mathematical concepts to study the complex Least Mean Square (LMS) algorithm, as used in adaptive array theory. 
\smallskip

In a previous work~\cite{nacht-adv}, inspired from the work of Brandwood \cite{Brandwood, WMcB}, the authors developed the bicomplex counterpart of the  C-R-Calculus and its application to bicomplex gradient operators and applied these new techniques to derive two bicomplex LMS algorithms extending both the classical real and complex LMS algorithms of Widrow. Here we discuss applications of this new bicomplex C-R-calculus (BC-R-Calculus) to the definitions and properties of bicomplex Hermite polynomials. \smallskip

The structure of our paper is as follows: in Section~\ref{complex} we remind the reader of the basic concepts and definitions of Complex-Real (C-R) Calculus. Then we review different notions of the classical complex case including polyanalytic functions, complex-Hermite polynomials, which are important examples of polyanalytic functions. We then briefly discuss the theory of complex gradient operators and we also explain several connections between polyanalytic functions and  C-R-regular functions. Section~\ref{bicomplex} collects different definitions and notations which will be useful in the rest of the paper. Here we review different notions related to the bicomplex algebra, the hyperbolic-valued modulus and bicomplex differential operators. In Section~\ref{BCR-Calc} we recall the bicomplex counterpart of the complex  C-R-Calculus and recall basic properties of classes of functions associated with the BC-R Calculus which was developed by the authors in~\cite{nacht-adv}. In Section~\ref{bhp} we present the important examples and properties of bicomplex Hermite polynomials and their relationships with the theory of polyanalytic functions.
%%%%%%%%%%%%%%%%%%%%%%%%%%%%%%%%%%%%%%%%%%%%%%%%%%%%%%%%%%%%
%%%%%%%%%%%%%%%%%%%%%%%%%%%%%%%%%%%%%%%%%%%%%%%%%%%%%%%%%%%%

\section{The complex case}
\label{complex}

This section is an expository reminder of several concepts from complex analysis, including polyanalytic function theory, It\^o or complex Hermite polynomials, C-R Calculus (or Wirtinger calculus).
We use the regular setting of a complex valued function $f:\Omega\subset \mathbb{C}\longrightarrow \mathbb{C}$, where $\Omega\subset \mathbb{C}$ is a domain. We consider $f$ to be real analytic with respect to its real and imaginary parts, therefore the operators $\displaystyle \frac{\partial}{\partial {z}}$ and $\displaystyle \frac{\partial}{\partial \overline{z}}$ are well-defined by

$$\displaystyle \frac{\partial}{\partial z}=\frac{1}{2}\left(\frac{\partial}{\partial x}-i\frac{\partial}{\partial y}\right), \quad\displaystyle \frac{\partial}{\partial \overline{z}}=\frac{1}{2}\left(\frac{\partial}{\partial x}+i\frac{\partial}{\partial y}\right);\quad z=x+iy.$$

Real differentiable functions in the kernel of $\displaystyle \frac{\partial}{\partial \overline{z}}$ are called {\em complex analytic, or, equivalently holomorphic} on their domains.

%%%%%%%%%%%%%%%%%%%%%%%%%%%%%%%%%%%%%

%%%%%%%%%%%%%%%%%%%%%%%%%%%%%%%%%%
\subsection{C-R Calculus in the complex case}
\label{Cplx CR}

We now introduce the concept of {\bf C-R analytic (or C-R regular)} functions which are the basis of study for {\it Complex-real or Wirtinger Calculus} i.e. C-R calculus for short~\cite{Brandwood, Kreutz}.
\begin{definition}
\label{CR-analytic}
Let $\Omega\subset \mathbb{C}$ be a complex domain, symmetric with respect to the real axis. We  say that a function $f:\Omega\longrightarrow \mathbb{C}$ is {\bf C-R analytic (or C-R regular) on $\Omega$} if there exists a complex analytic function of two complex variables $g:\Omega\times \Omega \longrightarrow \mathbb{C}, (z_1,z_2)\mapsto g(z_1,z_2)$ such that 
\begin{itemize}
\item[i)] $f(z)=g(z,\overline{z}), \quad \forall z\in \Omega$ ,
\\
\item[ii)] $\displaystyle \frac{\partial f}{ \partial z}=\left(\frac{\partial g}{\partial z_1} \right)_{z_1=z, \,z_2=\overline{z}}$ ,
\\
\item[iii)]  $\displaystyle \frac{\partial f}{ \partial \overline{z}}=\left(\frac{\partial g}{\partial z_2} \right)_{z_1=z, \, z_2=\overline{z}}.$
\end{itemize}
\end{definition}

\begin{example}
The function $\displaystyle f_1(z)=\Im(z)=\frac{z-\overline{z}}{2i}$ is a real valued C-R analytic function with $\displaystyle g_1(z_1,z_2)=\frac{z_1-z_2}{2i}.$ 
Also, the function $f_2(z)=z\overline{z}^2$ is a complex valued C-R analytic function with $g_2(z_1,z_2)=z_1z_2^2$. Note that neither $f_1$ nor $f_2$ are complex analytic.\\
Other interesting examples include $|z|^2$ and $\ln(z\bar{z})$.
\end{example}
We now write an example of derivatives of such functions.

\begin{example}
The following formulas hold for the derivatives of $|z|^{2k}$ and $|z|^{2k+1}$:
$$\frac{\partial |z|^{2k}}{ \partial \overline{z}}=kz|z|^{2k-2}; \quad \frac{\partial |z|^{2k+1}}{ \partial \overline{z}}=\frac{2k+1}{2}z|z|^{2k-1}, \quad k\ge 0 .$$
\end{example}

\begin{remark}
If $g$ is entire we have
\begin{equation*}\label{entire-cplx}
g(z_1,z_2)=\displaystyle \sum_{j=0}^{\infty}\sum_{k=0}^{\infty}a_{k,j} z_1^jz_2^k , \qquad \text{yielding} \qquad f(z)=\displaystyle \sum_{j=0}^{\infty}\sum_{k=0}^{\infty}a_{k,j} z^j\overline{z}^k ,
\end{equation*}
and the following forms for the derivatives follow:
\begin{equation*}\label{entire-cplx-derivatives}
 \frac{\partial f}{ \partial z}= \sum_{j=1}^{\infty}\left(\sum_{k=0}^{\infty}a_{k,j} \overline{z}^k \right) j \, z^{j-1} ,\qquad \text{and} \qquad 
 \frac{\partial f}{ \partial \overline{z}}=  \sum_{k=1}^{\infty}\left(\sum_{j=0}^{\infty}a_{k,j} z^k \right) k \, \overline{z}^{k-1} .
\end{equation*}
\end{remark}

In classic complex analysis we see that any function $f:\mathbb{C} \longrightarrow \mathbb{C}$ can be seen as a function from $\mathbb{R}\times\mathbb{R}$ to  $\mathbb{C}$ as $f(z)=f(x+ iy)=f(x,y)$. We re-formulate results of Brandwood~\cite{Brandwood}, for greater readability.

\begin{theorem} 
\label{Brand}
Let $f$ be C-R regular on $\mathbb{C}$, i.e., as in Definition~\ref{CR-analytic}, there exists a complex analytic function of two complex variables $g:\mathbb{C} \times \mathbb{C} \longrightarrow \mathbb{C}, (z_1,z_2)\mapsto g(z_1,z_2)$ such that 
$$f(z)=f(x,y)=g(z,\overline{z}); \quad z=x+iy.$$ Then, the differential operators follow the expected rules:  
\begin{eqnarray*} 
\frac{\partial g}{\partial z}&=\frac{1}{2}\left(\frac{\partial}{\partial x}-i\frac{\partial}{\partial y}\right)f(x,y),\\
 \frac{\partial g}{\partial \overline{z}}&=\frac{1}{2}\left(\frac{\partial}{\partial x}+i\frac{\partial}{\partial y}\right)f(x,y),
 \end{eqnarray*}

where $\displaystyle \frac{\partial g}{\partial z}:=\left(\frac{\partial g}{\partial z_1} \right)_{z_1=z, \,z_2=\overline{z}}, \qquad \displaystyle \frac{\partial g}{\partial \overline{z}}:=\left(\frac{\partial g}{\partial z_2} \right)_{z_1=z, \,z_2=\overline{z}}$.
\end{theorem}

The importance of C-R calculus has been demonstrated throughout the last century, with applications from function theory to neural networks.
%%%%%%%%%%%%%%%%%%%%%%%%%%%%%%%%%%
%%%%%%%%%%%%%%%%%%%%%%%%%%%%%%%%%%
\subsection{Classical polyanalytic functions}
\label{C-poly}
Under the conditions above we have:

\begin{definition}
\label{cplx-poly}
Let $\Omega\subset \C$ be a complex domain and $\displaystyle f\in\mathcal{C}^n(\Omega),$ i.e. $f$ real analytic of order $n$. If $f$ is in the kernel of a power $n\geq 1$ of the classical Cauchy-Riemann operator $\displaystyle \frac{\partial}{\partial \overline{z}}$, that is $$\displaystyle \frac{\partial^n}{\partial \overline{z}^n}f(z)=0, \quad \forall  z\in\Omega,$$
then $f$ is called a  {\bf polyanalytic function of order $n$} on $\Omega$.
The space of all  complex polyanalytic function of order $n$ on $\Omega$ is denoted by $ H_n(\Omega)$.
 \end{definition}
\smallskip
An interesting fact regarding these functions is that any polyanalytic function of order $n$ can be decomposed in terms of $n$ analytic functions so that we have a decomposition of the following form (see~\cite{Balk1991, balk_ency})
\begin{equation}
f(z)=\displaystyle \sum_{k=0}^{n-1}\overline{z}^kf_k(z),
\end{equation}
for which all $f_k$ are complex analytic functions on $\Omega$. In particular, expanding each analytic component using the series expansion theorem leads to an expression of the form

\begin{equation}\label{exp1}
f(z)=\displaystyle \sum_{k=0}^{n-1}\sum_{j=0}^{\infty}\overline{z}^kz^ja_{k,j},
\end{equation}
where $(a_{k,j})$ are complex coefficients.

In this paper we are also interested in the case where the expansion \eqref{exp1} is {\em of infinite order}, (this case was also considered in~\cite{IEOT_ACDS}),  we now write a new definition:
\begin{definition}
A function of the form
\begin{equation}\label{exp2}
f(z)=\displaystyle \sum_{k=0}^{\infty}\sum_{j=0}^{\infty}\overline{z}^kz^ja_{k,j},
\end{equation}
where the coefficients $a_{k,j}$ are non-zero for an infinite number of indices $k$, will be called a polyanalytic function of  {\bf infinite order}. 
\end{definition}

We note that such functions were discussed in \cite{Balk1991, balk_ency} as well, where they were mentioned as {\em conjugate analytic functions}.

\begin{proposition}
Any polyanalytic function (of fininte or infinte order) on a domain $\Omega$ symmetric with respect to the $x-$axis is C-R analytic. The converse is also true.
\end{proposition}

\begin{proof}
Let $f:\mathbb{C}\longrightarrow \mathbb{C}$ be a polyanalytic function of order $n=1,2,...$. We know by poly-decomposition that there exists unique analytic functions $f_0,....,f_{n-1}$ such that 
$$f(z)=\displaystyle \sum_{k=0}^{n-1}\overline{z}^kf_k(z); \quad \forall z\in \mathbb{C}.$$
Then, we consider the analytic function of two complex variables defined by $$g(z_1,z_2)=\displaystyle \sum_{k=0}^{n-1}z_2^kf_k(z_1), \quad \forall (z_1,z_2)\in \mathbb{C}^2.$$ 
It is clear that $g(z,\overline{z})=f(z)$ and then we can easily check that $f$ is C-R analytic.

The infinite order case is treated similarly and we leave the proof to the reader.\\
The converse easily follows from the definition.
\end{proof}
\begin{example}
We consider the function $f(z)=e^{|z|^2}$. It is clear that $f$ is C-R analytic with $$g(z,w):=\displaystyle \sum_{n=0}^\infty \frac{z^n w^n}{n!}=e^{zw};\quad \forall z,w\in \mathbb{C}.$$
However, $f$ is not polyanalytic of a finite order.
\end{example}

For $n=1,2,...$ we recall that polyanalytic Fock spaces of order $n$ can be defined as follows
$$\mathcal{F}_n(\mathbb{C}):=\left\lbrace g\in H_n(\mathbb{C}), \quad \frac{1}{\pi}\int_{\mathbb{C}}|g(z)|^2e^{-|z|^2}dA(z)<\infty \right\rbrace.$$
The reproducing kernel associated to the space $\mathcal{F}_n(\mathbb{C})$ is given by
\begin{equation}\label{Kn}
K_n(z,w)=e^{z\overline{w}}\displaystyle \sum_{k=0}^{n-1}\frac{(-1)^k}{k!}{n \choose k+1}|z-w|^{2k},
\end{equation}
for every $z,w\in\mathbb{C}.$\\

The poly-Bergman space $A^2_n(B(0,1))$ of polyanalytic functions of order $n$ in the complex unit disc $B(0,1)$ is given by $$\displaystyle A^2_n(B(0,1))=\lbrace{f\in H_n(B(0,1));\int_{B(0,1)}\vert{f(z)}\vert^2\, d A(z)<\infty}\rbrace.$$
We note that 
$A^2_n(B(0,1))$ is a reproducing kernel Hilbert space whose reproducing kernel is given by
\begin{equation}\label{BergKernel}
\displaystyle B_n(z,w)=\frac{n}{\pi(1-\overline{w}z)^{2n}}\sum_{k=0}^{n-1}(-1)^k{n \choose k+1}{n+k \choose n}\vert{1-\overline{w}z}\vert^{2(n-1-k)}\vert{z-w}\vert^{2k},
\end{equation}
for every $z,w\in \, B(0,1)$.

%%%%%%%%%%%%%%%%%%%%%%%%%%%%%%%%%%%%%%%%%%%%%%%%

\subsection{Complex Hermite Polynomials}

Important examples of polyanalytic functions were introduced by It\^o in~\cite{Ito}; these are called {\it complex Hermite polynomials} and they are used in the study of stochastic processes. It\^o also discovered their relation with complex multiple Wiener integral and we recall the definition of Hermite polynomials as follows:

\begin{definition}
Let $m,n$ be non-negative integers, then a polynomial of the form: 
\begin{equation}
H_{m,n}(z,\overline{z}):=\displaystyle \sum_{k=0}^{\min{(m,n)}}(-1)^k k! {m\choose k}{n \choose k}z^{m-k}\overline{z}^{n-k}
\end{equation}
is called a {\bf complex Hermite polynomial}.
\end{definition}

\begin{remark}
Alternatively, these polynomials can be defined via the Rodrigues formula:
\begin{equation}
H_{m,n}(z,\overline{z}):=(-1)^{m+n}e^{|z|^2}\frac{\partial^{m+n}}{\partial \overline{z}^m\partial z^n}(e^{-|z|^2})\,; \quad m,n\ge 0 .
\end{equation}
\end{remark}

The complex Hermite polynomials are very useful, as they provide an orthogonal basis for the $L^2$ space on $\mathbb{C}$ with respect to the classical Gaussian measure. The following integral formula holds
$$\displaystyle \int_{\mathbb{C}}H_{m,n}(z,\overline{z})\overline{H_{m',n'}(z,\overline{z})}e^{-|z|^2}dA(z)=\pi m!n! \delta_{((m,n);(m',n'))},$$
where $dA(z)$ is the usual Lebesgue measure on $\mathbb{C}$.

In~\cite{Ito} It\^o, as well as Ismail and Simeonov in~\cite{Ismail}, note the existence of a generating function in this case:

\begin{proposition}
For every $u, v\in \mathbb{C}$ the generating function given by the Hermite polynomials is
\begin{equation}
\displaystyle \sum_{m,n=0}^{\infty} H_{m,n}(z,\overline{z}) \frac{u^m v^n}{m! n!}=e^{uz+v\overline{z}-uv}, \quad \forall z\in \mathbb{C}.
\end{equation}

In particular, taking $u=w$ and $v=\overline{w}$ we obtain 
\begin{equation}
\displaystyle \sum_{m,n=0}^{\infty} H_{m,n}(z,\overline{z}) \frac{w^m \overline{w}^n}{m! n!}=e^{wz+\overline{w}\overline{z}-|w|^2}, \quad \forall z\in \mathbb{C}.
\end{equation}

\end{proposition}

\begin{remark}
The complex Hermite polynomials satisfy the following Appell-type properties:
\begin{equation}
\displaystyle \frac{\partial}{\partial z} H_{m,n}(z,\overline{z})=mH_{m-1,n}(z,\overline{z}), \quad \frac{\partial}{\partial \overline{z}} H_{m,n}(z,\overline{z})=nH_{m,n-1}(z,\overline{z}).
\end{equation} 
\end{remark}

Another important property of these polynomials is given by introducing the following two differential operators which are adjoints of each other with respect to the Gaussian $L^2$-norm
\begin{equation}
\displaystyle \mathcal{A}:=\frac{\partial}{\partial \overline{z}},\quad \mathcal{A}^*:=-\frac{\partial}{\partial z}+\overline{z}.
\end{equation}
These operators lead to a factorization of the Landau operator as in Proposition 7.1 in~\cite{Shige}.
\begin{definition}
The  {\bf Landau operator} is defined as:
\begin{equation}
\mathcal{G}:= -\frac{\partial^2}{\partial \overline{z} \partial z}+\overline{z}\frac{\partial}{\partial \overline{z}},
\end{equation} 
which is factorized by
$$\mathcal{G}=\mathcal{A}^*\mathcal{A}.$$
\end{definition}

Moreover, it turns out that complex Hermite polynomials are also eigenfunctions of the Landau operator $\mathcal{G}$,  and the following formula holds
\begin{equation}
\mathcal{G}(H_{m,n})(z,\overline{z})=nH_{m,n}(z,\overline{z}).
\end{equation}
Finally it is important to note that It\^o had already proven that the normalized complex Hermite polynomials $\psi_{m,n}=\frac{H_{m,n}}{\sqrt{m!n!}}$ form an orthonormal basis of $L^2(\mathbb{C},\frac{e^{-|z|^2}}{\pi}dA(z)).$

%%%%%%%%%%%%%%%%%%%%%%%%%%%%%%%%%%%%%%
%%%%%%%%%%%%%%%%%%%%%%%%%%%%%%%%%%%%%%

\section{Introduction to Bicomplex Numbers}
\label{bicomplex}

The algebra of bicomplex numbers was first introduced by Segre in~\cite{Segre}. During the past decades, a few isolated works analyzed either the properties of bicomplex numbers, or the properties of holomorphic functions defined on bicomplex numbers, and, without pretense of completeness, we direct the attention of the reader first to the to book of Price,~\cite{Price}, where a full foundation of the theory of multicomplex numbers was given, then to some of the works describing some analytic properties of functions in the field~\cite{alss, CSVV, DSVV, mltcplx}. Applications of bicomplex (and other hypercomplex) numbers can be also found in the works of Alfsmann, Sangwine, Gl\"{o}cker, and Ell~\cite{AG, AGES}.\\

We now introduce, in the same fashion as in~\cite{CSVV,bcbook,Price}, the key definitions and results for the case of bicomplex holomorphic functions of bicomplex variables. The algebra of bicomplex numbers is generated by two commuting imaginary units $\i$
and $\j$ and we will denote the bicomplex space by $\BC$.  The product of the two commuting units  $\i$ and $\j$ is denoted by $ \k := \i\j$ and we note that $\k$ is a hyperbolic unit, i.e. it is a unit which squares to $1$.  Because of these various units in $\BC$, there are several
different conjugations that can be defined naturally. We will make use of appropriate conjugations in this paper, and we refer the reader to~\cite{bcbook,mltcplx} for more information on bicomplex and multicomplex analysis.

\medskip
%%%%%%%%%%%%%%%%%%%%%%%%%%%%
%%%%%%%%%%%%%%%%%%%%%%%%%%%%
\subsection{Properties of the bicomplex algebra}

The bicomplex space, $\BC$, is not a division algebra, and it has two distinguished zero
divisors, $\e_1$ and $\e_2$, which are idempotent, linearly independent
over the reals, and mutually annihilating with respect to the
bicomplex multiplication:
\begin{align*}
  \e_1&:=\frac{1+\k}{2}\,,\qquad \e_2:=\frac{1-\k}{2}\,,\\
  \e_1 \cdot \e_2 &= 0,\qquad
  \e_1^2=\e_1 , \qquad \e_2^2 =\e_2\,,\\
  \e_1 +\e_2 &=1, \qquad \e_1 -\e_2 = \k\,.
\end{align*}
Just like $\{1,\mathbf{j} \},$ these form a basis of the complex algebra
$\BC$, which is called the {\em idempotent basis}. If we define the
following complex variables in $\C(\i)$:
\begin{align*}
  \beta_1 := z_1-\i z_2,\qquad \beta_2 := z_1+\i z_2\,,
\end{align*}
the $\C(\i)$--{\em idempotent representation} for $Z=z_1+\j z_2$ is
given by
\begin{align*}
  Z &= \beta_1\e_1+\beta_2\e_2\,.
\end{align*}

The $\C(\i)$--idempotent is the only representation for which
multiplication can be taken component-wise, as shown in the next lemma.

\begin{remark}
  \label{prop:idempotent}
  The addition and multiplication of bicomplex numbers can be realized
  component-wise in the idempotent representation above. Specifically,
  if $Z= a_1\,\e_2 + a_2\,\e_2$ and $W= b_1\,\e_1 + b_2\,\e_2 $ are two
  bicomplex numbers, where $a_1,a_2,b_1,b_2\in\C(\i)$, then
  \begin{eqnarray*}
    Z+W &=& (a_1+b_1)\,\e_1  + (a_2+b_2)\,\e_2   ,  \\
    Z\cdot W &=& (a_1b_1)\,\e_1  + (a_2b_2)\,\e_2   ,  \\
    Z^n &=& a_1^n \,\e_1  + a_2^n \,\e_2  .
  \end{eqnarray*}
  Moreover, the inverse of a bicomplex number 
  $Z=a_1\e_1 + a_2\e_2 $ is defined when $a_1 \cdot a_2 \neq 0$ and is given
  by
  $$
  Z^{-1}= a_1^{-1}\e_1 + a_2^{-1}\,\e_2 ,
  $$
  where $a_1^{-1}$ and $a_2^{-1}$ are the complex multiplicative
  inverses of $a_1$ and $a_2$, respectively.
\end{remark}

One can see this also by computing directly which product on the
bicomplex numbers of the form
\begin{align*}
  x_1 + \i x_2 + \j x_3 + \k x_4,\qquad x_1,x_2,x_3,x_4\in\R
\end{align*}
is component wise, and one finds that the only one with this property
is given by the mapping:
\begin{align}
  \label{shakira}
  x_1 + \i x_2 + \j x_3 + \k x_4 \mapsto ((x_1 + x_4) + \i (x_2-x_3), (x_1-x_4) + \i (x_2+x_3))\,,
\end{align}
which corresponds exactly with the idempotent decomposition
\begin{align*}
  Z = z_1 + \j z_2 = (z_1- \i z_2)\e_1 + (z_1+ \i z_2)\e_2\,,
\end{align*}
where $z_1 = x_1+ \i x_2$ and $z_2 = x_3+ \i x_4$.

\begin{remark}
These split the bicomplex space in $\BC=\mathbb C \mathbf{e}_1\bigoplus \mathbb C \mathbf{e}_2$, as:
\begin{equation}
  Z=z_1+\j z_2=(z_1-\i z_2)\mathbf{e}_1+(z_1+\i z_2)\mathbf{e}_2=\lambda_1\e_1+\lambda_2\e_2.
\end{equation}
\end{remark}

Simple algebra yields:
\begin{equation}
\begin{split}
  z_1&=\frac{\lambda_1+\lambda_2}{2}\, ,\\
    z_2&=\frac{\i(\lambda_1-\lambda_2)}{2}.
  \end{split}
  \end{equation}

Because of these various units in $\BC$, there are several
different conjugations that can be defined naturally and we will now define the conjugates in the bicomplex setting, as in~\cite{CSVV,bcbook}

\begin{definition} 
  \label{conj}
  For any $Z\in \BC$ we have the following three {\bf conjugates}:
  \begin{eqnarray*}
    \overline{Z}=\overline{z_1}+\j\overline{z_2}\, ,\\
     Z^{\dagger}=z_1-\j z_2\, ,\\
      Z^*=\overline{Z^{\dagger}}=\overline{z_1}-\j\overline{z_2}.
  \end{eqnarray*}
\end{definition}

\begin{remark}
\label{id-conj}
Moreover, following Definition~\ref{conj}, if we write $Z=\lambda_1 \mathbf{e}_1+\lambda_2\mathbf{e}_2$ in the idempotent representation,we have 
\begin{align*}
  Z &= \lambda_1 \e_1 + \lambda_2 \e_2  \, ,\\
 Z^\ast  &=  \overline{\lambda_1} \e_1 + \overline{\lambda_2} \e_2 \, ,\\
   Z ^\dagger&= \lambda_2 \e_1 + \lambda_1 \e_2  \, ,\\
   \overline{Z} &=  \overline{\lambda_2} \e_1 + \overline{\lambda_1} \e_2\,.
\end{align*}
\end{remark}
We refer the reader to~\cite{bcbook} for more details.
\bigskip

%%%%%%%%%%%%%%%%%%%%%%%%%%%%%%%%%%%%%%%%%%

%%%%%%%%%%%%%%%%%%%%%%%%%%%%%%%%%%%%%%%%%%
\subsection{Hyperbolic subalgebra and the hyperbolic-valued modulus}

A special subalgebra of $\BC$ is the set of hyperbolic numbers, denoted by $\D$.  This
algebra and the analysis of hyperbolic numbers have been studied, for
example, in~\cite{alss,bcbook,sobczyk} and we summarize
below only the notions relevant for our results.  A {\em hyperbolic number}
can be defined independently of $\BC$, by $\mathfrak{z}=x + \k y$,
with $x,y,\in\R$, $\k\not\in\R, \k^2=1$, and we denote by $\D$ the
algebra of hyperbolic numbers with the usual component--wise addition
and multiplication.  The hyperbolic {\em conjugate} of $\mathfrak z$
is defined by $\mathfrak{z}^\diamond := x - \k y$, and note that:
\begin{equation}
  \mathfrak{z}\cdot \mathfrak{z}^\diamond=x^2-y^2\in\R\,,
\end{equation}
which yields the notion of the square of the {\em modulus} of a
hyperbolic number $\mathfrak{z}$, defined by
$ |\mathfrak{z}|_{\D}^2:=\mathfrak{z}\cdot \mathfrak{z}^\diamond$.
\begin{remark}
It is worth noting that both $\overline{Z}$ and $Z^{\dagger}$ reduce to  $\mathfrak{z}^\diamond$ when $Z= \mathfrak{z}$. In particular $\e_2=\e_1^\diamond=\e_1^*=\e_1^{\dagger}$.
\end{remark}
Similar to the bicomplex case, hyperbolic numbers have a unique
idempotent representation with real coefficients:
\begin{align}
  \label{D_idempotent}
  \mathfrak{z}=s \e_1 + t \e_2 \,,
\end{align}
where, just as in the bicomplex case, $\displaystyle \e_1 = \frac{1}{2} \left( 1 + \k \right) $,
$\displaystyle \e_2= \frac{1}{2} \left( 1 - \k \right)$, and
$s:=x+y$ and $t:=x-y$. Note that $\e_1^\diamond=\e_2$ if we consider
$\D$ as a subset of $\BC$, as briefly explained in the remark above. We also observe that
$$
|\mathfrak{z}|_{\D}^2 = x^2 - y^2 = (x+y)(x-y) = st.
$$

The hyperbolic algebra $\D$ is a subalgebra of the bicomplex numbers
$\BC$ (see~\cite{bcbook} for details). Actually $\BC$ is the
algebraic closure of $\D$, and it can also be seen as the
complexification of $\D$ by using either of the imaginary unit $\i$ or
the unit $\j$.
\begin{definition}
Define the set $\D^+$ of {\bf non--negative hyperbolic numbers} by:
\begin{align*}
  \D^+ &= \left\{ x + \k y \, \big| \, x^2 - y^2 \geq 0,  x \geq 0 \right\}
       = \left\{ x + \k y \, \big| \, x \geq 0,  | y | \leq x \right\} \\
       &= \{ s \e_1 + t \e_2 \, \big| \, s, t \geq 0 \}.
\end{align*}
\end{definition}
\begin{remark}
As studied extensively in~\cite{alss}, one can define a partial order
relation defined on $\D$ by:
\begin{align}
  \label{po}
  \mathfrak{z}_1 \preceq \mathfrak{z}_2\qquad\text{if and only if}\qquad
  \mathfrak{z}_2-\mathfrak{z}_1\in\D^+,
\end{align}
and we will use this partial order to study the
{\em hyperbolic--valued} norm, which was first introduced and studied
in~\cite{alss}.
\end{remark}

The Euclidean norm $\|Z\|$ on $\BC$, when it is seen as
$\C^2(\i), \C^2(\j)$ or $\R^4$ is:
\begin{align*}
  \|Z\| = \sqrt{ | z_1 | ^2 + | z_2 |^2 \, } = \sqrt{ \Re\left( | Z |_\k^2 \right) \, } = \sqrt{
  \, x_1^2 + y_1^2 + x_2^2 + y_2^2 \, }.
\end{align*}
As studied in detail in~\cite{bcbook}, in idempotent
coordinates $Z=\lambda_1\e_1+\lambda_2\e_2$, the Euclidean norm becomes:
\begin{align}
  \label{Euclidean_idempotent}
  \|Z\| = \frac{1}{\sqrt2}\sqrt{|\lambda_1|^2 + |\lambda_2|^2}.
\end{align}

% For two bicomplex numbers $Z=z_1+\j z_2=x_1+x_2\i +x_3\j +x_4\k$ and
% $W=w_1+\j w_2=y_1+y_2\i +y_3\j +y_4\k$ one can define the inner
% product associated with the Euclidean norm to be
% $\ip{Z,W}_\R=\ip{z_1,w_1}_\R+\ip{z_1,w_1}_\R$, where
% $\ip{z_1,w_1}_\R=x_1y_1+x_2y_2$ is the usual inner product on $\R^2$.

It is easy to prove that
\begin{align}
  \|Z \cdot W\|  \leq  \sqrt{2} \left(\|Z\| \cdot  \|W\| \right),
\end{align}
and we note that this inequality is sharp since if $Z = W= \e_1$, one
has:
\begin{align*}
  \|\e_1 \cdot \e_1\| = \|\e_1\| = \frac{1}{\sqrt{ 2}} =  \sqrt{2}\, \|\e_1\| \cdot \|\e_1\|,
\end{align*}
and similarly for $\e_2$.

\begin{definition}
One can define a
{\bf hyperbolic-valued} norm for $Z=z_1+\j z_2 = \lambda_1 \e_1+\lambda_2\e_2$
by:
\begin{align*}
  \| Z\|_{\D_+} := |\lambda_1|\e_1 + |\lambda_1|\e_2 \in\D^+.
\end{align*}
\end{definition}
It is shown in~\cite{alss} that this definition obeys the corresponding properties
of a norm, i.e. $ \| Z\|_{\D_+}=0$ if and only if $ Z=0$,
it is multiplicative, and it
respects the triangle inequality with respect to the order introduced
above. 

%%%%%%%%%%%%%% 
%%%%%%%%%%%%%% 
\subsection{Finsler--type Norm}

Another real valued norm that can be useful in this setting is the one
found by multiplying all conjugates of a bicomplex number.  This norm has been first considered in~\cite{vm} and let us recall that it is a fourth order Finsler-type norm defined by
$$
|Z|_{\mathcal F}^4 := Z\overline{Z} Z^\ast Z^\dagger.
$$
For $Z=z_1+\j z_2 = \varphi_1 + \i\varphi_2$, with $z_1,z_2\in\C_\i$
and $\varphi_1,\varphi_2\in\C_\j$, recall that:
\begin{align*}
  \overline{Z} &:= \overline{z_1} + \j\overline{z_2} = \varphi_1 - \i \varphi_2,\qquad
            Z^\dagger := z_1 - \j z_2 = \overline{\varphi_1} - \i\overline{\varphi_2},\\ 
  Z^\ast &:= \left(\overline{Z}\right)^\dagger =
           \overline{\left(Z^\dagger\right)} = \overline{z_1} - \j\overline{z_2}
           = \overline{\varphi_1} + \i\overline{\varphi_2}
           = \mathfrak{z}_1 - \i\, \mathfrak{z}_2
           = \mathfrak{w}_1 - \j\, \mathfrak{w}_2
\end{align*}
so 
\begin{align*}
  |Z|_{\mathcal F}^4 &= Z\overline{Z} Z^\ast Z^\dagger 
               = (z_1 + \j z_2)(\overline{z_1} + \j\overline{z_2})(\overline{z_1} - \j\overline{z_2})(z_1 - \j z_2)
               = (z_1^2 + z_2^2)(\overline{z_1}^2 + \overline{z_2}^2).
\end{align*}
The corresponding (square of the) bicomplex moduli are defined as:
\begin{align*}
  |Z|_{\i}^2 &:= Z \cdot Z^\dagger = z_1^2 + z_2^2 \in \C_{\i}\\ 
  |Z|_{\j}^2 &:= Z \cdot \overline{Z} = \varphi_1^2  + \varphi_2^2 \in \C_{\j}\\
  |Z|_{\k}^2 &:= Z \cdot Z^\ast = \mathfrak z_1^2  + \mathfrak z_2^2  
               =  \mathfrak w_1^2   + \mathfrak w_2^2 \in\D.
\end{align*}
so
\begin{align*}
  |Z|_{\mathcal F}^4 &= |Z|_{\i}^2\cdot|\overline{Z}|_{\i}^2 = \left||Z|_\j^2\right|_\i^2= ||Z|_\i^2|_\j^2 = ||Z|_\k^2|_\i^2.
\end{align*}

In this case we will obtain a modulus of order $4$ as $\mathcal F(Z)$ is a
real number.  This number is positive when $Z$ is not a divisor of $0$
and it is equal to $0$ when $Z\in\mathfrak{S}$, i.e. it is a complex
multiple of $\e$ or $\edag$.

From Remark~\ref{id-conj} we have:
$$\overline{Z}=\overline{\lambda_2} \mathbf{e}_1+\overline{\lambda_1}\mathbf{e}_2; \quad Z^*=\overline{\lambda_1} \mathbf{e}_1+\overline{\lambda_2}\mathbf{e}_2; \quad Z^{\dagger}=\lambda_2 \mathbf{e}_1+\lambda_1 \mathbf{e}_2.$$

which brings the following form for the three bicomplex moduli of order $2$ and the modulus of order $4$ to :
\begin{lemma}
In the idempotent representation we have:
$$|Z|_{\i}^2= \lambda_1\lambda_2;\quad |Z|_{\j}^2=\lambda_1\overline{\lambda_2} \mathbf{e}_1+\overline{\lambda_1}\lambda_2\mathbf{e}_2; \quad |Z|_{\k}^2=|\lambda_1|^2\mathbf{e}_1+|\lambda_2|^2\mathbf{e}_2 ;\quad \mathcal F(Z)=|\lambda_1|^2 \cdot |\lambda_2|^2. $$
\end{lemma}

For proofs and more details about these moduli one can see~\cite{vm} and~\cite{SVV}.

%%%%%%%%%%%%%%%%%%%%%%%%%%%%%%%%%%%%%%%%%%%%%%%%%%%%%%%

%%%%%%%%%%%%%%%%%%%%%%%%%%%%%%%%%%
\section{A review of the Bicomplex-Real (BC-R) calculus}
\label{BCR-Calc}

In this section we reiterate the extension of the study of C-R (or Wirtinger) Calculus to the bicomplex setting the authors~\cite{nacht-adv}.

\subsection{Bicomplex differential operators}

We first define the bicomplex differential operators arising from the rich underlying structure of the space.

In the case of the bicomplex variable, $Z=z_1+\mathbf{j}z_2$, we have $$z_1=\frac{Z+Z^\dagger}{2}, \quad z_2=\frac{\mathbf{j}}{2}(Z^\dagger-Z),$$
and 
$$\overline{z_1}=\frac{\overline{Z}+Z^*}{2}, \quad \overline{z_2}=\frac{\mathbf{j}}{2}(Z^*-\overline{Z}).$$
On the other hand, for $z_1=x_1+\i x_2$ and $z_2=x_3+\i x_4$ we obtain

$$x_1=\frac{Z+Z^\dagger+Z^*+\overline{Z}}{4}, \quad x_2=\frac{\mathbf{i}}{4}(Z^*+\overline{Z}-Z-Z^\dagger),$$
and $$x_3=\frac{\mathbf{j}}{4}(Z^*+Z^{\dagger}-\overline{Z}-Z), \quad x_4=\frac{\mathbf{k}}{4}(Z+Z^*-\overline{Z}-Z^\dagger).$$
The bicomplex differential operators with respect to the various bicomplex conjugates (see \cite{alss, CSVV}) are  $$\partial_{Z}:= \partial_{z_1}-\mathbf{j} \partial_{z_2}=\partial_{x_1}-\mathbf{i}\partial{x_2}-\mathbf{j}\partial_{x_3}+\mathbf{k}\partial{x_4}=\partial_{\lambda_1}\mathbf{e}_1+\partial_{\lambda_2}\mathbf{e}_2\, ,$$
 $$\partial_{\overline{Z}}:= \partial_{\overline{z_1}}-\mathbf{j} \partial_{\overline{z_2}}=\partial_{x_1}+\mathbf{i}\partial{x_2}-\mathbf{j}\partial_{x_3}-\mathbf{k}\partial{x_4}= \partial_{\overline{\lambda_2}}\mathbf{e}_1+\partial_{\overline{\lambda_1}}\mathbf{e}_2 \, ,$$
  $$\partial_{Z^*}:= \partial_{\overline{z_1}}+\mathbf{j} \partial_{\overline{z_2}}=\partial_{x_1}+\mathbf{i}\partial{x_2}+\mathbf{j}\partial_{x_3}+\mathbf{k}\partial{x_4}=\partial_{\overline{\lambda_1}}\mathbf{e}_1+\partial_{\overline{\lambda_2}}\mathbf{e}_2 \, ,$$
  $$\partial_{Z^\dagger}:=\partial_{z_1}+\mathbf{j} \partial_{z_2}=\partial_{x_1}-\mathbf{i}\partial{x_2}+\mathbf{j}\partial_{x_3}-\mathbf{k}\partial{x_4}=\partial_{\lambda_2}\mathbf{e}_1+\partial_{\lambda_1}\mathbf{e}_2 \, .$$

 We recall  the definition of bicomplex analyticity,  a more recent concept in hypercomplex analysis, more details in~\cite{alv1, alv2, alss, bcbook, Price}.
 
 \begin{definition}
 Let $\Omega$ be a domain in $\BC$, we say that a function $F:\Omega\longrightarrow \BC$ is bicomplex holomorphic if and only if $F$ is in the kernel of the last three differential operators described above, i.e.  $\partial_{\overline{Z}}, \partial_{Z^*},  \partial_{Z^\dagger}$.
 \end{definition}
 
 We invite the reader to the following references~\cite{alss, CSVV, bcbook, Price} for more details. Following these works, we can see that a bicomplex holomorphic function will admit a convergent  power series representation at each point in $\Omega$. We point out that in Section 7.6 from~\cite{bcbook}, we can see that a function is bicomplex analytic if and only if it can be split in the idempotent representation in two functions each depending on a single idempotent variable. 
 
 As in~\cite{CSVV}, the concept of bicomplex holomorphy can be extended to several bicomplex variables and we have:
 \begin{definition}
 Let $\Omega$ be a domain in $\BC^n$, we say that a function $F:\Omega\longrightarrow \BC$ is a bicomplex holomorphic function of several bicomplex variables if and only if $F$ is in the kernel of the family of the $3n$ differential operators (three for each variable)  $\partial_{\overline{Z_l}}, \partial_{Z_l^*},  \partial_{Z_l^\dagger}$, for any $1\le l \le n$, where the $n$ bicomplex vector variable is denoted by $Z=(Z_1,\dots, Z_n)$.
 \end{definition}

We will use this last definition throughout our paper and for more details on the theory of functions of several bicomplex variables we refer the reader to~\cite{CSVV} .

\subsection{Notions of BC-R calculus}

In this subsection we define the counterpart of C-R-Calculus (or Wirtinger) in the bicomplex space.
We start by introducing the notion of BC-R analytic functions and discuss some of their main properties. In what follows we use the notion of open bicomplex domains, namely domains of the form $\Omega_1\mathbf{e}_1+\Omega_2\mathbf{e}_2$ where $\Omega_1,\, \Omega_2$ are complex domains.

\begin{definition}
\label{B-C-R-Omega}
Let $\Omega \,\subset\,\BC$ be a bicomplex domain, symmetric with respect to all conjugations, i.e. if $Z\in\Omega$ then $\overline{Z},Z^*,Z^{\dagger}\in\Omega$.
We  say that a function $\displaystyle f:\Omega\longrightarrow \mathbb{BC}$ is {\bf BC-R analytic (or BC-R regular)} on $\Omega$ if there exists a $\BC$-analytic function of four bicomplex variables $$g:\Omega^4 \longrightarrow \mathbb{BC}, (Z_1,Z_2,Z_3,Z_4)\mapsto g(Z_1,Z_2,Z_3,Z_4)$$ such that 
\begin{itemize}
\item[i)] $f(Z)=g(Z,\overline{Z},Z^*,Z^\dagger), \quad \forall Z\in \Omega$,
\\
\item[ii)] $\displaystyle \frac{\partial f}{ \partial Z}=\left(\frac{\partial g}{\partial Z_1} \right)_{Z_1=Z, \,Z_2=\overline{Z},\,Z_3=Z^*\,Z_4=Z^\dagger}$ ,
\\
\item[iii)]  $\displaystyle \frac{\partial f}{ \partial \overline{Z}}=\left(\frac{\partial g}{\partial Z_2} \right)_{Z_1=Z, \,Z_2=\overline{Z},\,Z_3=Z^*\,Z_4=Z^\dagger}$,
\item[iv)] $\displaystyle \frac{\partial f}{ \partial Z^*}=\left(\frac{\partial g}{\partial Z_3} \right)_{Z_1=Z, \,Z_2=\overline{Z},\,Z_3=Z^*\,Z_4=Z^\dagger}$,
\item[v)] $\displaystyle \frac{\partial f}{ \partial Z^\dagger}=\left(\frac{\partial g}{\partial Z_4} \right)_{Z_1=Z, \,Z_2=\overline{Z},\,Z_3=Z^*\,Z_4=Z^\dagger}$.
\end{itemize}
\end{definition}
\begin{remark} Just as the Complex-real calculus generalizes notions of analytic complex functions, the bicomplex-real calculus generalizes notions of bicomplex analyticity. We point out that in Section 7.6 from~\cite{bcbook}, we can see that a function is bicomplex analytic if and only if it can be split in the idempotent representation in two functions each depending on a single idempotent variable. 
This type of analyticity generalizes this concept.
\end{remark}
\begin{example}
The Finsler-type norm defined by $$f(Z)=|Z|_{\mathcal F}^4 := Z\bar{Z} Z^\ast Z^\dagger$$ is an entire B-C-R function on $\mathbb{BC}$, with $$ g(Z_1,Z_2,Z_3,Z_4)=Z_1Z_2Z_3Z_4.$$
We have $$\displaystyle \frac{\partial f}{ \partial Z}=\bar{Z} Z^\ast Z^\dagger, \quad  \frac{\partial f}{ \partial \overline{Z}}=Z Z^\ast Z^\dagger,$$
and $$\frac{\partial f}{ \partial Z^*}=Z\bar{Z}  Z^\dagger, \quad \frac{\partial f}{ \partial Z^\dagger}=Z\bar{Z} Z^\ast .$$
\end{example}
\begin{proposition}
Let $f,h:\mathbb{BC}\longrightarrow \mathbb{BC}$ be two BC-R analytic functions and $\lambda \in\mathbb{BC} $. Then, the sum $f+h$ and multiplication $f\lambda$ are also BC-R analytic.
\end{proposition}
\begin{proof}
Follows standard arguments.
\end{proof}
\begin{remark}
For $\Omega=\BC$ in Definition~\ref{B-C-R-Omega}, we call such functions {\em BC-R entire} and the set of all BC-R entire functions is a vector space over $\mathbb{BC}$ which is denoted by $\mathcal{H}_{CR}(\mathbb{BC})$.
\end{remark}
We note that bicomplex polyanalytic functions of finite order were considered in \cite{Ghanmi2022} using finite sums involving different bicomplex conjugates. 

\begin{definition}
\label{bcplx-poly-finite}
A bicomplex valued function $f:\mathbb{BC}\longrightarrow \mathbb{BC}$ which satisfies
$$\displaystyle \frac{\partial^l f}{\partial \overline{Z}^l}=\frac{\partial^k f}{\partial (Z^*)^k}=\frac{\partial^q f}{\partial (Z^\dagger)^q}=0,$$
where $l,k$ and $q$ are strictly positive integers, is called a {\bf bicomplex polyanalytic function of finite multi-order $(l, k, q)$}.
\end{definition}
\begin{remark}
For example, in~\cite{Ghanmi2022}, Proposition 3.8 shows that a function is as in Definition~\ref{bcplx-poly-finite} if and only if
$$g(Z):=\displaystyle \sum_{m=0}^{l-1} \sum_{n=0}^{k-1}\sum_{p=0}^{q-1} g_{m,n,p}(Z)\, {\overline{Z}}^m(Z^*)^n(Z^{\dagger})^p,$$
where $g_{m,n,p}$ are bicomplex holomorphic functions. It is also worth noting that bicomplex analytic functions are a special case of bicomplex polyanalytic ones.
\end{remark}

 We introduce the class of global bicomplex polyanalytic functions of infinite order as follows
\begin{definition}
A bicomplex valued function $f:\mathbb{BC}\longrightarrow \mathbb{BC}$ is called a {\bf global bicomplex polyanalytic of infinite order} if it can be represented as a power series with respect to the variables $Z, \overline{Z},Z^*, Z^{\dagger}$ so that we have 
\begin{equation}
\displaystyle f(Z)=\sum_{m,n,p,q=0}^{\infty} a_{m,n,p,q} \, Z^m{\overline{Z}}^n(Z^*)^p(Z^{\dagger})^q ;
\end{equation}
where $(a_{m,n,p,q})_{m,n,p,q\geq 0}$ are suitable bicomplex coefficients, non-zero for an infinite number of indices $n,p$ or $q$. 
\end{definition}

In~\cite{nacht-adv} we have proven the following theorem on the relationship between BC- analyticity and bicomplex polyanalycity.
\begin{theorem}
A bicomplex valued function $f:\mathbb{BC}\longrightarrow \mathbb{BC}$ is global polyanalytic of finite or infinite order if and only if it is a BC-R analytic function.
\end{theorem}
These global polyanalytic types of functions can take many forms, depending on whether they are in the kernels of the operators 
$\displaystyle \frac{\partial}{\partial \overline{Z}}, \, \frac{\partial}{\partial Z^*},\, \frac{\partial^q f}{\partial Z^\dagger}$ or their powers, yielding power series in terms of specific conjugates, or exhibiting finite order polyanalytic type behavior, respectively.

When we consider the extension of the complex Laplacian to the bicomplex setting, there exists more than one form, depending on the choice of bicomplex conjugate. For example one can define three of these Laplacians as follows

\begin{definition}
We define the Laplacians with respect to $\i, \j$ and $\k$ as follows
\begin{equation}
\displaystyle \Delta_\i=\frac{\partial^2}{\partial Z \partial \overline{Z}};
\end{equation}

\begin{equation}
\displaystyle \Delta_\j:=\frac{\partial^2}{\partial Z \partial Z^{\dagger}};
\end{equation}

\begin{equation}
\displaystyle \Delta_\k=\frac{\partial^2}{\partial Z \partial Z^*}.
\end{equation}

\end{definition}
 The following example relates to the Laplacian given by $Z^*$ and in~\cite{nacht-adv} we have proven the following result concerning the logarithm of a bicomplex function $f(Z)$, denoted by $\mathsf{Ln}(f(Z))$. This result extends a complex one proven in  Exercise $4.2.23$ from~\cite{daniel_problem_book}, also found in~\cite{evgrafov}, to the bicomplex setting.
\begin{proposition}
Let $h$ be a $\mathbb{BC}$ holomorphic function, we consider the function defined by $$f(Z)=||h(Z)||_{\mathbf{k}}^{2}g(Z).$$ We denote by $\Delta_{\k}$ the bicomplex Lalplacian given by $$\Delta_{\k}:=\displaystyle \frac{\partial^2}{\partial Z\partial Z^*}.$$
Then, it holds that \begin{equation}
\Delta_{\k} \mathsf{Ln} f(Z)=\Delta_{\k} \mathsf{Ln} g(Z).
\end{equation}
\end{proposition}

\subsection{A type of bicomplex Fock space}

An useful type of polyanalyticity is the one defined by the $*-$conjugation and we have:
\begin{definition}
\label{bcplx-*-poly} 
If $f$ is in the kernel of $\displaystyle \frac{\partial}{ \partial \overline{Z}}$,  $\displaystyle \frac{\partial}{ \partial Z^{\dag}}$ and in the kernel of a power $n\geq 1$ of the $*$ operator $\displaystyle \frac{\partial}{\partial Z^*}$, that is $$\displaystyle \frac{\partial^n}{\partial {Z^*}^n}f(Z)=0, \quad \forall  Z\in { \Omega}\subset\BC,$$
where $\Omega=\Omega_1\e_1+\Omega_2\e_2$ is a product (or split) domain in $\BC$,
then $f$ is called a bicomplex $*$ polyanalytic function of order $n$ on $\Omega$.\\
The space of all  bicomplex $*$ polyanalytic function of order $n$ is denoted by $ H_n^{split}(\Omega)$.
 \end{definition}

Just as in the case of complex polyanalytic function of order $n$ can be decomposed in terms of $n$ analytic functions so that we have a decomposition of the following form

\begin{theorem} Let $f$ be a bicomplex $*$ polyanalytic function of order $n$ on $\Omega\subset\BC$ a split domain, then we have:
\begin{equation}
f(Z)=\displaystyle \sum_{l=0}^{n-1}{Z^*}^lf_l(Z),
\end{equation}
for which all $f_l$ are bicomplex analytic functions on $\Omega$. In particular, expanding each analytic component using the series expansion theorem leads to an expression of this form.
\end{theorem}

In~\cite{nacht-adv} we extended the definition of a polyanalytic Fock space in Section~\ref{complex}, to the bicomplex case.

For $n=1,2,...$ we recall that polyanalytic Fock spaces of order $n$ can be defined as follows
$$\mathcal{F}_n(\mathbb{BC}):=\left \lbrace g\in H_n^{split}(\mathbb{BC}), \quad   g(\lambda_1 \e_1 +\lambda_2\e_2)=g_1(\lambda_1)\e_1+g_2(\lambda_2)\e_2, \, g_1,g_2\in \mathcal{F}_n(\mathbb{C})\right\rbrace,$$ and we have:
\begin{proposition}
The reproducing kernel associated to the space $\mathcal{F}_n(\mathbb{BC})$ is given by
\begin{equation}\label{BC_Kn}
\mathsf{K}_n(Z,W)=e^{Z{W^*}}\displaystyle \sum_{l=0}^{n-1}\frac{(-1)^l}{l!}{n \choose l+1}|Z-W|_{\k}^{2l},
\end{equation}
for every $Z,W\in\mathbb{BC}.$ 
\end{proposition}
\begin{proof} 
Please see~\cite{nacht-adv}.
\end{proof}

Another special case follows:
\begin{definition}
\label{bcplx-poly-Bergman}
The bicomplex $*$-poly-Bergman space $\mathcal{A}_n(\mathbb{K})$ of polyanalytic functions of order $n$ in the bicomplex product-type unit ball denoted here by $\mathbb{K}=B(0,1)\e_1+B(0,1)\e_2$ is given by $$\mathcal{A}_n(\mathbb{K}):=\left \lbrace g\in H_n^{split}(\mathbb{K}), \quad   g(\lambda_1 \e_1 +\lambda_2\e_2)=g_1(\lambda_1)\e_1+g_2(\lambda_2)\e_2, \, g_1,g_2\in A^2_n(B(0,1))\right\rbrace,$$
where

\end{definition}

\begin{remark}
we note also that 
$\mathcal{A}_n(\mathbb{K})$ is a reproducing kernel Hilbert space whose reproducing kernel is given by
\begin{equation}
\displaystyle \mathsf{B}_n(Z,W)=\frac{n}{\pi(1-W^*Z)^{2n}}\sum_{\ell=0}^{n-1}(-1)^\ell {n \choose \ell+1}{n+\ell \choose n}\vert{1-W^*Z}\vert_{\mathbf{k}}^{2(n-1-\ell)}\vert{Z-W}\vert^{2\ell}_{\mathbf{k}}
\end{equation}
for every $Z,W\in\mathbb{K}$.
\end{remark}

So, as in the case of the bicomplex poly-Fock space we can show that for $Z=\lambda_1\mathbf{e}_1+\lambda_2\mathbf{e}_2$ and $W=\mu_1\mathbf{e}_1+\mu_2\mathbf{e}_2$ we have  

\begin{equation}\label{idem_BC_BergKn}
\mathsf{B}_n(Z,W)=B_n(\lambda_1,\mu_1)\e_1+B_n(\lambda_2,\mu_2)\e_2,
\end{equation}
where $B_n(\lambda_1,\mu_1)$ and $B_n(\lambda_2,\mu_2)$ are given in~\eqref{BergKernel}.

\begin{remark}
In the case $n=1$ we recover the Bergman kernel of the bicomplex Bergman space, found in Corollary 6.2 of~\cite{Quirroga}. In this case the functions will be bicomplex analytic.
\end{remark}

%%%%%%%%%%%%%%%%%%%%%%%%%%%%%%%%%%%%%%%
%%%%%%%%%%%%%%%%%%%%%%%%%%%%%%%%%%%%%%%

\section{Bicomplex Hermite Polynomials}
\label{bhp}
\setcounter{equation}{0}
Bicomplex Hermite polynomials are important families of BC-R analytic functions. We will show that they generate spaces of bicomplex functions with certain properties.

\subsection{Bicomplex Hermite polynomials of the first kind}
\begin{definition}
Let $m,n=0,1,....$ we define three bicomplex Hermite polynomials of the first kind (or It\^o bicomplex polynomials) with respect to each conjugate
\begin{enumerate}
\item The $\overline{Z}$-bicomplex Hermite polynomials are given by $$
H_{m,n}(Z,\overline{Z}):=\displaystyle \sum_{k=0}^{\min{(m,n)}}(-1)^k k! {m\choose k}{n \choose k}Z^{m-k}\overline{Z}^{n-k};
$$
\item The $Z^*$-bicomplex Hermite polynomials are given by$$
H_{m,n}(Z,Z^*):=\displaystyle \sum_{k=0}^{\min{(m,n)}}(-1)^k k! {m\choose k}{n \choose k}Z^{m-k}(Z^*)^{n-k};
$$
\item The $Z^{\dagger}$-bicomplex Hermite polynomials are given by $$
H_{m,n}(Z,Z^{\dagger}):=\displaystyle \sum_{k=0}^{\min{(m,n)}}(-1)^k k! {m\choose k}{n \choose k}Z^{m-k}(Z^{\dagger})^{n-k}.
$$
\end{enumerate}

\end{definition}
\begin{remark}
For every $m,n=0,1,...$ we have 
$$H_{m,0}(Z,\overline{Z})=H_{m,0}(Z,Z^*)=H_{m,0}(Z,Z^\dagger)=Z^m;$$
 and $$H_{0,n}(Z,\overline{Z})=\overline{Z}^n, \quad H_{0,n}(Z,Z^*)=(Z^*)^n, \quad H_{0,n}(Z,Z^\dagger)=(Z^\dagger)^n.$$
\end{remark}
\begin{proposition}
The bicomplex Hermite polynomials of the first kind $H_{m,n}(Z,\overline{Z}), H_{m,n}(Z,Z^*)$ and $H_{m,n}(Z,Z^{\dagger})$ are B-C-R entire for every $m,n=0,1, \cdots$.
\end{proposition}
\begin{proof}
Follows standard arguments.
\end{proof}
Now we shall state different results on these bicomplex polynomials of the first kind. We start by the following technical result 
\begin{proposition}
\label{Hermite_split}
For every $Z=z_1+z_2\mathbf{j}=\lambda_1\mathbf{e}_1+\lambda_2\mathbf{e}_2\in \mathbb{BC}$ we have 
\begin{equation}
H_{m,n}(Z,Z^*)=H_{m,n}(\lambda_1,\overline{\lambda_1})\mathbf{e}_1+H_{m,n}(\lambda_2,\overline{\lambda_2})\mathbf{e}_2
\end{equation}
\end{proposition}
\begin{proof}
We note that using the expression of classical complex Hermite polynomials with respect to the variables $(\lambda_1,\overline{\lambda_1})$ and $(\lambda_2,\overline{\lambda_2})$ we have 
$$H_{m,n}(\lambda_1,\overline{\lambda_1})=\displaystyle \sum_{k=0}^{\min{(m,n)}}(-1)^k k! {m\choose k}{n \choose k}\lambda_{1}^{m-k}\overline{\lambda_1}^{n-k}, $$
and $$ H_{m,n}(\lambda_2,\overline{\lambda_2})=\displaystyle \sum_{k=0}^{\min{(m,n)}}(-1)^k k! {m\choose k}{n \choose k}\lambda_{2}^{m-k}\overline{\lambda_2}^{n-k}.$$
Hence, it follows 
\begin{align*}
\displaystyle H_{m,n}(\lambda_1,\overline{\lambda_1})\mathbf{e}_1+H_{m,n}(\lambda_2,\overline{\lambda_2})\mathbf{e}_2&=\sum_{k=0}^{\min(m,n)}(-1)^k k! {m\choose k}{n \choose k}\left( \lambda_{1}^{m-k}\overline{\lambda_1}^{n-k}\mathbf{e_1}+\lambda_{2}^{m-k}\overline{\lambda_2}^{n-k}\mathbf{e}_2\right)\\
&=\sum_{k=0}^{\min(m,n)}(-1)^k k! {m\choose k}{n \choose k}\left( \lambda_{1}^{m-k}\mathbf{e_1}+\lambda_{2}^{m-k}\mathbf{e}_2\right)\left(\overline{\lambda_1}^{n-k}\mathbf{e_1}+\overline{\lambda_2}^{n-k}\mathbf{e}_2\right)\\
&=\sum_{k=0}^{\min(m,n)}(-1)^k k! {m\choose k}{n \choose k}\left( \lambda_{1}\mathbf{e_1}+\lambda_{2}\mathbf{e}_2\right)^{m-k}\left(\overline{\lambda_1}\mathbf{e_1}+\overline{\lambda_2}\mathbf{e}_2\right)^{n-k}\\
&=\sum_{k=0}^{\min{(m,n)}}(-1)^k k! {m\choose k}{n \choose k}Z^{m-k}(Z^*)^{n-k}\\
&=H_{m,n}(Z,Z^*)\\
\end{align*}
This ends the proof.
\end{proof}

\begin{proposition}
The It\^o-bicomplex Hermite polynomials satisfy the following Appell properties
$$\displaystyle \frac{\partial}{\partial Z} H_{m,n}(Z,Z^*)=mH_{m-1,n}(Z,Z^*), \qquad \displaystyle \frac{\partial}{\partial Z^*} H_{m,n}(Z, Z^*)=nH_{m,n-1}(Z,Z^*).$$

\end{proposition}
\begin{proof}
We provide explicit computations with respect to $Z^*$-variable. The same reasoning works for the $Z$-variable. In fact, we note that applying the previous propositio we have $$H_{m,n}(Z,Z^*)=H_{m,n}(\lambda_1,\overline{\lambda_1})\mathbf{e}_1+H_{m,m}(\lambda_2,\overline{\lambda_2})\mathbf{e}_2.$$
Then we use the expression of the differential operator $\frac{\partial}{\partial Z^*}$ in terms of the idempotent representation to get 
\begin{align*}
\displaystyle \frac{\partial}{\partial Z^*} H_{m,n}(Z,Z^*)&=\frac{\partial}{\partial  \overline{\lambda_1}}H_{m,n}(\lambda_1,\overline{\lambda_1})\mathbf{e}_1+\frac{\partial}{\partial \overline{\lambda_2}}H_{m,n}(\lambda_2,\overline{\lambda_2})\mathbf{e}_2\\
&=nH_{m,n-1}(\lambda_1,\overline{\lambda_1})\mathbf{e}_1+nH_{m,n-1}(\lambda_2,\overline{\lambda_2})\\
&=nH_{m,n-1}(Z,Z^*).\\
\end{align*}
\end{proof}

\begin{theorem}[Rodrigues formula]
The following representation hold for the bicomplex Hermite polynomials of the first kind
\begin{equation}
H_{m,n}(Z,Z^*)=(-1)^{m+n}e^{|Z|^{2}_{\mathbf{k}}}\frac{\partial^{m+n}}{\partial (Z^*)^m\partial Z^n}(e^{-|Z|^{2}_{\mathbf{k}}}); \quad m,n=0,1,...
\end{equation}

\end{theorem}
\begin{proof}
Let $Z=z_1+z_2\mathbf{j}=\lambda_1\mathbf{e}_1+\lambda_2\mathbf{e}_2$.
We note that $$\displaystyle e^{-|Z|^{2}_{\mathbf{k}}}=e^{-|\lambda_1|^{2}}\mathbf{e}_1+e^{-|\lambda_2|^{2}} \mathbf{e}_2,$$
and $$\displaystyle \frac{\partial^{m+n}}{\partial (Z^*)^m\partial Z^n}=\frac{\partial^{m+n}}{\partial \overline{\lambda_1}^m \partial \lambda_1^n }\mathbf{e}_1+ \frac{\partial^{m+n}}{\partial \overline{\lambda_2}^m \partial \lambda_2^n }\mathbf{e}_2$$
So
$$\displaystyle \frac{\partial^{m+n}}{\partial (Z^*)^m\partial Z^n}(e^{-|Z|^{2}_{\mathbf{k}}})=\frac{\partial^{m+n}}{\partial \overline{\lambda_1}^m \partial \lambda_1^n }\left(e^{-|\lambda_1|^{2}}\right) \mathbf{e}_1+ \frac{\partial^{m+n}}{\partial \overline{\lambda_2}^m \partial \lambda_2^n }\left(e^{-|\lambda_2|^{2}}\right)\mathbf{e}_2.$$
 Then, applying Proposition~\ref{Hermite_split} combined with the classical Rodrigues formula for complex Hermite polynomials, we obtain: 
\begin{align*}
H_{m,n}(Z,Z^*)&=H_{m,n}(\lambda_1,\overline{\lambda_1})\mathbf{e}_1+H_{m,m}(\lambda_2,\overline{\lambda_2})\mathbf{e}_2,\\
&=(-1)^{m+n}e^{|\lambda_1|^{2}}\frac{\partial^{m+n}}{\partial \overline{\lambda_1}^m\partial \lambda_1^n}(e^{-|\lambda_1|^2})\mathbf{e}_1+(-1)^{m+n}e^{|\lambda_1|^{2}}\frac{\partial^{m+n}}{\partial \overline{\lambda_2}^m\partial \lambda_2^n}(e^{-|\lambda_2|^2})\mathbf{e}_2,\\
&=(-1)^{m+n}e^{|Z|^{2}_{\mathbf{k}}}\frac{\partial^{m+n}}{\partial (Z^*)^m\partial Z^n}(e^{-|Z|_{\mathbf{k}}^2}). \\
\end{align*}

\end{proof}

\begin{theorem}
The Hermite polynomials of first kind in $\overline{Z}$ and $Z^{\dagger}$ also satisfy Rodrigues's formula as follows:
\begin{equation}
H_{m,n}(Z,\overline{Z})=(-1)^{m+n}e^{|Z|^{2}_{\mathbf{j}}}\frac{\partial^{m+n}}{\partial \overline{Z}^m\partial Z^n}(e^{-|Z|_{\mathbf{j}}^2}); \quad m,n\ge 0 ,
\end{equation}
and
\begin{equation}
H_{m,n}(Z,Z^{\dagger})=(-1)^{m+n}e^{|Z|^{2}_{\mathbf{i}}}\frac{\partial^{m+n}}{\partial (Z^\dagger)^m\partial Z^n}(e^{-|Z|_{\mathbf{i}}^2}); \quad m,n\ge 0 .
\end{equation}

\end{theorem}
\begin{proof}
The proofs follow arguments similar with the classical complex proof, using the Leibniz rule and an inductive process.
\end{proof}

\begin{theorem}[Generating functions]\label{GF*}
For every $U,V, Z\in \mathbb{BC}$ we have 
\begin{equation}
\displaystyle \sum_{m,n=0}^{\infty} H_{m,n}(Z,Z^*) \frac{U^m V^n}{m! n!}=e^{UZ+VZ^*-UV}.
\end{equation}
\end{theorem}
\begin{proof}
We write $Z=z_1+z_2\mathbf{j}=\lambda_1\mathbf{e}_1+\lambda_2\mathbf{e}_2$,  $U=u_1+u_2\mathbf{j}=u'_1\mathbf{e}_1+u'_2\mathbf{e}_2$ and $V=v_1+v_2\mathbf{j}=v'_1\mathbf{e}_1+v'_2\mathbf{e}_2$. 
Observe that using Proposition~\ref{Hermite_split} we obtain
\begin{align*}
\displaystyle  \sum_{m,n=0}^{\infty} H_{m,n}(Z,Z^*) \frac{U^m V^n}{m! n!}&= \sum_{m,n=0}^{\infty} \left( H_{m,n}(\lambda_1,\overline{\lambda_1})\mathbf{e}_1+H_{m,n}(\lambda_2,\overline{\lambda_2})\mathbf{e}_2\right) \frac{U^m V^m}{m! n!}\\
&=\left(\sum_{m,n=0}^{\infty} H_{m,n}(\lambda_1,\overline{\lambda_1}) \frac{(u'_{1})^{m}(v'_{1})^{n}}{m! n!}\right) \mathbf{e}_1+\left(\sum_{m,n=0}^{\infty} H_{m,n}(\lambda_2,\overline{\lambda_2}) \frac{(u'_{2})^{m}(v'_{2})^{n}}{m! n!}\right) \mathbf{e}_2\\
&=e^{u'_1\lambda_1+v'_1\overline{\lambda_1}-u'_1v'_1}\mathbf{e}_1+e^{u'_2\lambda_2+v'_2\overline{\lambda_2}-u'_2v'_2}\mathbf{e}_2\\
&=e^{UZ+VZ^*-UV}.\\
\end{align*}
\end{proof}
\begin{remark} Let $Z, U\in \mathbb{BC}$.
We observe that if $V=U^*$ in Theorem \ref{GF*}  we get a hyperbolic valued function given by
\begin{equation}
\displaystyle \sum_{m,n=0}^{\infty} H_{m,n}(Z,Z^*) \frac{U^m (U^*)^n}{m! n!}=e^{UZ+(UZ)^*-|U|^{2}_{\mathbf{k}}}.
\end{equation}
\end{remark}

Replacing the $*$ bicomplex conjugate by any of the other two specific to this theory (namely $\overline{\,}$ and $\dagger$ ) similar argument can be made about the other first order Hermite polynomials we consider, however the proof is more aligned with the proof found in the original papers by  It\^o~\cite{Ito}. The derivation follows along similar lines and we have:

\begin{lemma}
We can consider also the following generating functions
\begin{equation}\label{I1}
\displaystyle \sum_{m,n=0}^{\infty} H_{m,n}(Z,\overline{Z}) \frac{U^m V^n}{m! n!}=e^{UZ+V\overline{Z}-UV}, \quad \forall Z\in \mathbb{BC};
\end{equation} and 
\begin{equation}\label{I2}
\displaystyle \sum_{m,n=0}^{\infty} H_{m,n}(Z,Z^\dagger) \frac{U^m V^n}{m! n!}=e^{UZ+V Z^\dagger-UV}, \quad \forall Z\in \mathbb{BC}.
\end{equation}
\end{lemma}
\begin{proof}
First we observe that the bicomplex exponential function satisfies the following functional equation $$e^{Z+W}=e^Z e^W, \quad \forall Z,W \in \mathbb{BC}.$$ 
As a consequence, we note that the calculations of the bicomplex formulas \eqref{I1} and \eqref{I2} follow a similar argument used by It\^o in Theorem 12 of his original paper \cite{Ito}. 
The remaining details are left to the reader.
\end{proof}
%{\color{blue} For this lemma, let us explore the case when $V$ is one of the appropriate conjugates of $U$}

\begin{definition}
One can define the full bicomplex $L^2$ space, denoted $L_k^2(\mathbb{BC}, e^{-Z Z^*})$ to be the space of functions $f$.... such that $f(Z)=f(\lambda_1\mathbf{e}_1+\lambda_2 \mathbf{e}_2)=f_1(\lambda_1,\lambda_2) \mathbf{e}_1 + f_2(\lambda_1, \lambda_2) \mathbf{e}_2,$ where $f_1$ and $f_2$ are $L_2$ integrable in $\mathbb{C}^2$ with Gaussian measure 
$e^{-(|\lambda_1|^2+ |\lambda_2|^2)}$.
\end{definition}

\begin{definition}
One can define the split-bicomplex $L^2$ space, denoted $L_{split}^2(\mathbb{BC}, e^{-Z Z^*})$ to be the space of split functions $f$, i.e. such that $f(Z)=f(\lambda_1\mathbf{e}_1+\lambda_2 \mathbf{e}_2)=f_1(\lambda_1) \mathbf{e}_1 + f_2(\lambda_2) \mathbf{e}_2,$ where $f_1$ and $f_2$ are $L^2$ integrable in $\mathbb{C}$ with Gaussian measure $e^{-|\lambda |^2}$.
\end{definition}

\begin{theorem}
The bicomplex Hermite polynomials $\lbrace{H_{m,n}(Z,Z^*);\quad m,n=0,1,\cdots}\rbrace$ form an orthogonal basis of the split-bicomplex Hilbert space $L_{split}^2(\mathbb{BC}, e^{-Z Z^*})$. Moreover, we have 
$$||H_{m,n}||_{L_{split}^2(\mathbb{BC}, e^{-Z Z^*})}^{2}=m!n!$$
\end{theorem}
\begin{proof}
Let $Z=z_1+\mathbf{j}z_2=\lambda_1\mathbf{e}_1+\lambda_2\mathbf{e}_2$ such that $$f(Z)=f_1(\lambda_1)\mathbf{e}_1+f_2(\lambda_2)\mathbf{e}_2$$
where $f_1, f_2$ are in $L^2(\mathbb{C},e^{-|w|^2}dA(w))$. Since complex Hermite polynomials form an orthogonal basis of $L^2(\mathbb{C},e^{-|w|^2}dA(w))$ we have 
\begin{align*}
f(Z)&=\left(\sum_{m,n=0}^{\infty} H_{m,n}(\lambda_1,\overline{\lambda_1})\alpha_{m,n} \right)\mathbf{e}_1+\left( \sum_{m,n=0}^{\infty} H_{m,n}(\lambda_2,\overline{\lambda_2})\beta_{m,n}\right)\mathbf{e}_2\\
&=\sum_{m,n=0}^{\infty}( H_{m,n}(\lambda_1,\overline{\lambda_1})\mathbf{e}_1+ H_{m,n}(\lambda_2,\overline{\lambda_2})\mathbf{e}_2)(\alpha_{m,n}\mathbf{e}_1+\beta_{m,n}\mathbf{e}_2) \\
&=\sum_{m,n=0}^{\infty}H_{m,n}(Z,Z^*)\gamma_{m,n}\,,\\
\end{align*}
with $\gamma_{m,n}=\alpha_{m,n}\mathbf{e}_1+\beta_{m,n}\mathbf{e}_2$ are coefficients in $\mathbb{BC}$.
\end{proof}

\begin{definition}[Bicomplex Landau operators]

We can define three Bicomplex Landau operators with respect to the three Laplacians $\displaystyle \Delta_\i,  \Delta_\j$ and $\Delta_\k$.
Here we define the bicomplex Landau operator given by 
\begin{equation}
\displaystyle \mathcal{G}_*:=-\Delta_\k+Z^*\frac{\partial}{\partial Z^*}.
\end{equation}
\end{definition}
\begin{lemma}
Just as in the complex case, it follows that this Landau operator admits a factorization
\begin{equation}
\displaystyle \mathsf{A}:=\frac{\partial}{\partial Z^*},\quad \mathsf{B}:=-\frac{\partial}{\partial Z}+Z^*.
\end{equation}
\end{lemma}
\begin{proof}
This is an immediate consequence of the definition of the BC-Landau operator.
\end{proof}

\begin{theorem}[Eigenfunctions of the BC-Landau operators]
For every $m,n=0,1,\cdots$ the bicomplex Hermite polynomials $H_{m,n}(Z,Z^*)$ are eigenfunctions of the bicomplex Landau operator $\mathcal{G}_*$. In other words, we have

\begin{equation}
\mathcal{G}_*(H_{m,n})(Z,Z^*)=nH_{m,n}(Z,Z^*).
\end{equation}

\end{theorem}
\begin{proof}
We note that 
\begin{align*}
\displaystyle \mathcal{G}_*:&= -\Delta_\k+Z^*\frac{\partial}{\partial Z^*} \\
&= \left(-\frac{\partial^2}{\partial \lambda_1 \partial \overline{\lambda_1}}+\overline{\lambda_1}\frac{\partial}{\partial \overline{\lambda_1}}\right)\mathbf{e}_1+\left(-\frac{\partial^2}{\partial \lambda_2 \partial \overline{\lambda_2}}+\overline{\lambda_2}\frac{\partial}{\partial \overline{\lambda_2}}\right)\mathbf{e}_2\\
\end{align*}
Moreover, we know by the classical complex case that 
$$\displaystyle  \left(-\frac{\partial^2}{\partial \lambda_1 \partial \overline{\lambda_1}}+\overline{\lambda_1}\frac{\partial}{\partial \overline{\lambda_1}}\right)H_{m,n}(\lambda_1,\overline{\lambda_1})=nH_{m,n}(\lambda_1,\overline{\lambda_1}), $$
and $$\displaystyle  \left(-\frac{\partial^2}{\partial \lambda_2 \partial \overline{\lambda_2}}+\overline{\lambda_2}\frac{\partial}{\partial \overline{\lambda_2}}\right)H_{m,n}(\lambda_2,\overline{\lambda_2})=nH_{m,n}(\lambda_2,\overline{\lambda_2}). $$
Hence, by splitting Hermite polynomials using Proposition \ref{Hermite_split} we obtain 
\begin{align*}\displaystyle 
\mathcal{G}_*(H_{m,n})(Z,Z^*)&=\left(-\frac{\partial^2}{\partial \lambda_1 \partial \overline{\lambda_1}}+\overline{\lambda_1}\frac{\partial}{\partial \overline{\lambda_1}}\right)H_{m,n}(\lambda_1,\overline{\lambda_1})\mathbf{e}_1+ \left(-\frac{\partial^2}{\partial \lambda_2 \partial \overline{\lambda_2}}+\overline{\lambda_2}\frac{\partial}{\partial \overline{\lambda_2}}\right)H_{m,n}(\lambda_2,\overline{\lambda_2}) \mathbf{e}_2\\
&=nH_{m,n}(\lambda_1,\overline{\lambda_1})\mathbf{e}_1+ nH_{m,n}(\lambda_2,\overline{\lambda_2})\mathbf{e}_2\\
&=nH_{m,n}(Z,Z^*).\\
\end{align*}
\end{proof}
\subsection{Bicomplex Hermite polynomials of the second kind}
In this section we consider bicomplex Hermite polynomials of the second kind. These polynomials can be related to the bicomplex field norm as follows:
\begin{definition}
Let $m,n, p, q\in \mathbb{N}$, the bicomplex Hermite polynomials of the second kind are defined by the following expression
%$$
%H_{m,n,p,q}(Z,\overline{Z},Z^*,Z^{\dagger}):=\displaystyle \sum_{k=0}^{\min{(m,n,p,q)}}(-1)^k k! {m\choose k}{n \choose k}{p\choose k}{m\choose k}Z^{m-k}\overline{Z}^{n-k}(Z^*)^{p-k}(Z^{\dagger})^{q-k}.
%$$
\begin{align*}
\displaystyle H_{m,n,p,q}(Z,\overline{Z},Z^*,Z^{\dagger})&:=(-1)^q m!n!p!q! Z^{q-n}(Z^*)^{q-p}(\overline{Z})^{q-m} \sum_{\ell=0}^{p} \sum_{k=0}^{m} \sum_{s=0}^{n} (-1)^{\ell+k+s}{p \choose \ell} {m \choose k} {n \choose s}  \\
&\times {q+\ell+k+s \choose n } \frac{||Z||^{4(k+\ell+s)}_{\mathcal{F}}}{(q+\ell+k+s)!}.\\
\end{align*}
\end{definition}
\begin{remark}
The bicomplex Hermite polynomials of the second kind are also B-C-R entire.
\end{remark}
\begin{theorem}[Rodrigues formula]
The following representation hold for the bicomplex Hermite polynomials of the second kind

$$H_{m,n,p,q}(Z,\overline{Z},Z^*,Z^{\dagger})= (-1)^{m+n+p+q}e^{|Z|_{\mathcal F}^4}\frac{\partial^{m+n+p+q}}{\partial \overline{Z}^m\partial Z^n \partial (Z^*)^p\partial (Z^\dagger)^q}(e^{-|Z|_{\mathcal F}^4}); \quad m,n,p,q=0,1,...$$
\end{theorem}
\begin{proof}
We observe that 
$$\frac{\partial^q}{\partial (Z^\dagger)^q}(e^{-|Z|_{\mathcal F}^4})=(-1)^q(Z\overline{Z}Z^*)^q e^{-|Z|_{\mathcal F}^4}$$
Setting $\psi_q(Z)=(Z\overline{Z}Z^*)^q$ we obtain 

$$\displaystyle \frac{\partial^{p+q}}{\partial(Z^*)^p\partial(Z^\dagger)^q}(e^{-|Z|_{\mathcal F}^4})=(-1)^q \sum_{\ell=0}^{p}{p \choose \ell} \frac{\partial^{p-\ell}}{\partial (Z^*)^{p-\ell}} \psi_q(Z)\frac{\partial^{\ell} }{\partial (Z^*)^\ell}(e^{-|Z|^4_\mathcal{F}})$$
Then, we develop the computations and get
\begin{align*}
\displaystyle \frac{\partial^{p-\ell}}{(Z^*)^{p-\ell}}\psi_q(Z)&=Z^q\overline{Z}^q \frac{\partial^{p-\ell}}{\partial (Z^*)^{p-\ell}}((Z^*)^q)\\
&=Z^q\overline{Z}^q \frac{q!}{(q-(p-\ell))!}(Z^*)^{q-(p-\ell)}\\
&=Z^q\overline{Z}^q \frac{q!}{(q-p+\ell)!}(Z^*)^{q-p+\ell}\\
\end{align*}
On the other hand we have also 
$$\displaystyle \frac{\partial^\ell}{\partial (Z^*)^\ell}\left(e^{-|Z|^{4}_{\mathcal{F}}}\right)=(-1)^\ell (Z\overline{Z}Z^\dagger)^\ell e^{-|Z|_{\mathcal{F}}^{4}}.$$
So, this leads to the following calculations 
$$\displaystyle \frac{\partial^{p+q}}{\partial(Z^*)^p\partial(Z^\dagger)^q}(e^{-|Z|_{\mathcal F}^4})=(-1)^q \sum_{\ell=0}^{p}{p \choose \ell} (-1)^\ell \frac{q!}{(q-p+\ell)!}Z^{q+\ell}\overline{Z}^{q+\ell}(Z^*)^{q-p+\ell} (Z^\dagger)^\ell e^{-|Z|_{\mathcal{F}}^{4}}$$
%At this stage we will apply the operator $\displaystyle \frac{\partial^m}{\partial \overline{Z}^m}$ and obtain the following
%\begin{align*}
%\displaystyle \frac{\partial^{m+p+q}}{\partial \overline{Z}^m\partial(Z^*)^p\partial(Z^\dagger)^q}(e^{-|Z|_{\mathcal F}^4})&=(-1)^q \sum_{\ell=0}^{p}{p \choose \ell} (-1)^\ell \frac{q!}{(q-p+\ell)!}Z^{q+\ell}(Z^*)^{q-p+\ell} (Z^\dagger)^\ell \frac{\partial^m}{\partial \overline{Z}^m}\left(\overline{Z}^{q+\ell} e^{-|Z|_{\mathcal{F}}^{4}} \right)\\
%\end{align*}
%Then, using similar arguments we have also 
%$$ \displaystyle \frac{\partial^m}{\partial \overline{Z}^m}\left(\overline{Z}^{q+\ell} e^{-|Z|_{\mathcal{F}}^{4}} \right)=\sum_{k=0}^{m}(-1)^k{m \choose k}\frac{(q+\ell)!}{(q+\ell-m+k)!}\overline{Z}^{q+\ell-m+k}(ZZ^*Z^\dagger)^ke^{-|Z|^{4}_{\mathcal{F}}}$$
Hence, developing the computations we obtain 

\begin{align*}
\displaystyle \frac{\partial^{m+p+q}}{\partial \overline{Z}^m\partial(Z^*)^p\partial(Z^\dagger)^q}(e^{-|Z|_{\mathcal F}^4})&=(-1)^q \sum_{\ell=0}^{p} \sum_{k=0}^{m} {p \choose \ell} {m \choose k} \frac{ (-1)^{\ell+k} q!(q+\ell)!}{(q-p+\ell)! (q+\ell-m+k)!}Z^{q+\ell}(Z^*)^{q-p+\ell} (Z^\dagger)^\ell \\
&\times \overline{Z}^{q+\ell-m+k}(ZZ^*Z^\dagger)^ke^{-|Z|^{4}_{\mathcal{F}}}\\
\end{align*}
Thus we have 

\begin{align*}
\displaystyle \frac{\partial^{m+p+q}}{\partial \overline{Z}^m\partial(Z^*)^p\partial(Z^\dagger)^q}(e^{-|Z|_{\mathcal F}^4})&=(-1)^q \sum_{\ell=0}^{p} \sum_{k=0}^{m} {p \choose \ell} {m \choose k} \frac{ (-1)^{\ell+k} q!(q+\ell)!}{(q-p+\ell)! (q+\ell-m+k)!}Z^{k+q+\ell}(Z^*)^{q-p+k+\ell}  \\
&\times (Z^\dagger)^{\ell+k} \overline{Z}^{q+\ell-m+k}e^{-|Z|^{4}_{\mathcal{F}}}\\
\end{align*}
Hence we obtain 
\begin{align*}
\displaystyle \frac{\partial^{m+n+p+q}}{\partial Z^n\partial \overline{Z}^m\partial(Z^*)^p\partial(Z^\dagger)^q}(e^{-|Z|_{\mathcal F}^4})&=(-1)^q \sum_{\ell=0}^{p} \sum_{k=0}^{m} {p \choose \ell} {m \choose k} \frac{ (-1)^{\ell+k} q!(q+\ell)!}{(q-p+\ell)! (q+\ell-m+k)!}(Z^*)^{q-p+k+\ell}  \\
&\times (Z^\dagger)^{\ell+k} \overline{Z}^{q+\ell-m+k} \frac{\partial^n}{\partial Z^n}\left( Z^{k+q+\ell}e^{-|Z|^{4}_{\mathcal{F}}}\right).\\
\end{align*}
Then, we calculate $\displaystyle \frac{\partial^n}{\partial Z^n}\left( Z^{k+q+\ell}e^{-|Z|^{4}_{\mathcal{F}}}\right)$ and obtain

\begin{align*}
\displaystyle \frac{\partial^{m+n+p+q}}{\partial Z^n\partial \overline{Z}^m\partial(Z^*)^p\partial(Z^\dagger)^q}(e^{-|Z|_{\mathcal F}^4})&=(-1)^q \sum_{\ell=0}^{p} \sum_{k=0}^{m} \sum_{s=0}^{n}{p \choose \ell} {m \choose k} {n \choose s} \frac{ (-1)^{\ell+k+s} q!(q+\ell)!(q+\ell+k)!}{(q-p+\ell)! (q+\ell+k-m)!}  \\
&\times \frac{1}{(q+\ell+k+s-n)!}  (Z^*)^{q-p+k+\ell+s}(Z^\dagger)^{\ell+k+s} \overline{Z}^{q+\ell+k+s-m} Z^{k+q+\ell+s-n}\\
\end{align*}

Hence, developing further the computations we get 
\begin{align*}
\displaystyle \frac{\partial^{m+n+p+q}}{\partial Z^n\partial \overline{Z}^m\partial(Z^*)^p\partial(Z^\dagger)^q}(e^{-|Z|_{\mathcal F}^4})&=(-1)^q m!n!p!q! Z^{q-n}(Z^*)^{q-p}(\overline{Z})^{q-m} \sum_{\ell=0}^{p} \sum_{k=0}^{m} \sum_{s=0}^{n}(-1)^{\ell+k+s}{p \choose \ell} {m \choose k} {n \choose s}  \\
&\times {q+\ell+k+s \choose n } \frac{||Z||^{4(k+\ell+s)}_{\mathcal{F}}}{(q+\ell+k+s)!}\\
\end{align*}
\end{proof}
\begin{remark}
An interesting $\mathcal{F}$-Bicomplex Laplace operator which can be considered in this setting is defined by
\begin{equation}
\displaystyle \Delta_\mathcal{F}:=\frac{\partial^4}{\partial Z \partial \overline{Z} \partial Z^* \partial Z^{\dagger}}.
\end{equation}
In forth-coming works we will study further properties of this operator and its relation with bicomplex Hermite polynomials of the second kind.
\end{remark}

%\begin{theorem}[Eigenfunctions of $\mathcal{F}$-Landau operator]

%\end{theorem}
%\begin{proof}

%\end{proof}
%\begin{theorem}[Generating functions]
 %We have 
 %\begin{equation}
%\displaystyle \sum_{m,n,p,q=0}^{\infty} H_{m,n, p,q}(Z,\overline{Z},Z^*,Z^\dagger) \frac{U^m V^nW^pY^q}{m! n!p!q!}=, \quad \forall Z\in \mathbb{BC};
%\end{equation}

%\end{theorem}
%\textbf{Question:} three variable Hermite polynomials ??

%%%%%%%%%%%%%%%%%%%%%%%%%%%%%%%%%%%%%%%%%%%%%%%%%%%%%%%%%%%%

\bibliographystyle{plain}

\end{document}